\documentclass[final,12pt]{elsarticle}
\usepackage{srcltx}
\usepackage{eurosym}
\usepackage{mathtools}
\usepackage{amsmath}
\usepackage{amsfonts}
\usepackage{amssymb}
\usepackage{amsthm}
\usepackage{graphicx}
\usepackage{mathrsfs}
\usepackage{xcolor}
\usepackage{exscale}
\usepackage{latexsym}

\usepackage[colorlinks,plainpages=true,pdfpagelabels,hypertexnames=true,colorlinks=true,pdfstartview=FitV,linkcolor=blue,citecolor=red,urlcolor=black]{hyperref}
\PassOptionsToPackage{unicode}{hyperref}
\PassOptionsToPackage{naturalnames}{hyperref}
\usepackage{enumerate}
\usepackage[shortlabels]{enumitem}
\usepackage{bookmark}
\usepackage{wasysym}
\usepackage{esint}
\usepackage[ddmmyyyy]{datetime}
\usepackage[margin=2cm]{geometry}
\makeatletter
\g@addto@macro\normalsize{%
	\setlength\abovedisplayskip{4pt}
	\setlength\belowdisplayskip{4pt}
	\setlength\abovedisplayshortskip{4pt}
	\setlength\belowdisplayshortskip{4pt}
}
\usepackage[capitalize,nameinlink]{cleveref}
\crefname{section}{Section}{Sections}
\crefname{subsection}{Subsection}{Subsections}
\crefname{condition}{Condition}{Conditions}
\crefname{hypothesis}{Hypothesis}{Conditions}
\crefname{assumption}{Assumption}{Assumptions}
\crefname{lemma}{Lemma}{Lemmas}
\crefname{definition}{Definition}{Definitions}

\crefformat{equation}{\textup{#2(#1)#3}}
\crefrangeformat{equation}{\textup{#3(#1)#4--#5(#2)#6}}
\crefmultiformat{equation}{\textup{#2(#1)#3}}{ and \textup{#2(#1)#3}}
{, \textup{#2(#1)#3}}{, and \textup{#2(#1)#3}}
\crefrangemultiformat{equation}{\textup{#3(#1)#4--#5(#2)#6}}%
{ and \textup{#3(#1)#4--#5(#2)#6}}{, \textup{#3(#1)#4--#5(#2)#6}}%
{, and \textup{#3(#1)#4--#5(#2)#6}}

\Crefformat{equation}{#2Equation~\textup{(#1)}#3}
\Crefrangeformat{equation}{Equations~\textup{#3(#1)#4--#5(#2)#6}}
\Crefmultiformat{equation}{Equations~\textup{#2(#1)#3}}{ and \textup{#2(#1)#3}}
{, \textup{#2(#1)#3}}{, and \textup{#2(#1)#3}}
\Crefrangemultiformat{equation}{Equations~\textup{#3(#1)#4--#5(#2)#6}}%
{ and \textup{#3(#1)#4--#5(#2)#6}}{, \textup{#3(#1)#4--#5(#2)#6}}%
{, and \textup{#3(#1)#4--#5(#2)#6}}

\crefdefaultlabelformat{#2\textup{#1}#3}

\newtheorem{theorem}{Theorem}

\newtheorem{lemma}[theorem]{Lemma}

\newtheorem{definition}[theorem]{Definition}

\def\CC{{\rm \kern.24em \vrule width.02em height1.4ex depth-.05ex \kern-.26emC}}

\def\TagOnRight

\def\AA{{it I} \hskip-3pt{\tt A}}

\def\QQ{\rlap {\raise 0.4ex \hbox{$\scriptscriptstyle |$}} {\hskip -0.1em Q}}

\makeatletter
\newcommand{\vo}{\vec{o}\@ifnextchar{^}{\,}{}}
\makeatother

\def\YYint#1#2#3{{\setbox0=\hbox{$#1{#2#3}{\iint}$}
		\vcenter{\hbox{$#2#3$}}\kern-.50\wd0}}

\def\XXint#1#2#3{{\setbox0=\hbox{$#1{#2#3}{\int}$}
		\vcenter{\hbox{$#2#3$}}\kern-.50\wd0}}

\makeatletter
\def\namedlabel#1#2{\begingroup
	\def\@currentlabel{#2}%
	\label{#1}\endgroup
}
\makeatother
\makeatletter
\newcommand{\rmh}[1]{\mathpalette{\raisem@th{#1}}}
\newcommand{\raisem@th}[3]{\hspace*{-1pt}\raisebox{#1}{$#2#3$}}
\makeatother

\newcounter{desccount}

\newcommand{\descref}[2]{\hyperref[#1]{\textnormal{\textcolor{black}{}\textcolor{blue}{\bf #2}\textcolor{black}{}}}}

\newcommand{\dref}[2]{\hyperref[#1]{\textcolor{black}{(}\textcolor{blue}{\bf #2}\textcolor{black}{)}}}
\newcommand{\be} {\begin{eqnarray}}
\newcommand{\ee} {\end{eqnarray}}
\newcommand{\Bea} {\begin{eqnarray*}}
	\newcommand{\Eea} {\end{eqnarray*}}
\newcommand{\pa} {\partial}
\newcommand{\re}{\mathbb{R}}

\newcommand{\al} {\alpha}
\newcommand{\rr}{\rightarrow}

\newcommand{\dip}{\displaystyle}

\newcommand{\e}  {\epsilon}
\newcommand{\De} {\Delta}

\newcommand{\f}{\infty}

\newcommand{\lab} {\label}

\newcommand{\vr} {\varrho}

\newcommand{\domm}{\int\limits_{\tau_1}^{\tau_2}\int\limits_{\Omega}} 
\newcommand{\V}{\textbf{v}_\e}
\newcommand{\dom}{\int\limits_{\Omega}}                     

\newcommand{\nb}{\nabla}

\DeclareMathOperator{\dv}{div}

\newcounter{whitney}
\refstepcounter{whitney}

\newcounter{ineqcounter}
\refstepcounter{ineqcounter}
\makeatletter
\def\ps@pprintTitle{%
	\let\@oddhead\@empty
	\let\@evenhead\@empty
	\def\@oddfoot{}%
	\let\@evenfoot\@oddfoot}
\makeatother
\usepackage[doublespacing]{setspace}
\usepackage[titletoc,toc,page]{appendix}

\makeatletter
\newcommand{\refcheckize}[1]{%
	\expandafter\let\csname @@\string#1\endcsname#1%
	\expandafter\DeclareRobustCommand\csname relax\string#1\endcsname[1]{%
		\csname @@\string#1\endcsname{##1}\wrtusdrf{##1}}%
	\expandafter\let\expandafter#1\csname relax\string#1\endcsname
}
\makeatother

\refcheckize{\cref}
\refcheckize{\Cref}


\makeatletter
\newcommand{\mainsectionstyle}{%
	\renewcommand{\@secnumfont}{\bfseries}
	\renewcommand\section{\@startsection{section}{2}%
		\z@{.5\linespacing\@plus.7\linespacing}{-.5em}%
		{\normalfont\bfseries}}%
}
\makeatother
\usepackage{pgf,tikz}
\usetikzlibrary{arrows}
\usetikzlibrary{decorations.pathreplacing}

\usepackage{xpatch}
\xpatchcmd{\MaketitleBox}{\hrule}{}{}{}
\xpatchcmd{\MaketitleBox}{\hrule}{}{}{}

\usepackage[utf8]{inputenc}
\date{}
\begin{document}
	\begin{frontmatter}
		
		\title{A sufficient condition for uniqueness of weak solutions of the incompressible Euler system	 \vspace{+2ex}
		} 
		\author[myaddress]{Shyam Sundar Ghoshal}
		\ead{ghoshal@tifrbng.res.in}

		\author[myaddress]{Animesh Jana}
		\ead{animesh@tifrbng.res.in }
		
		\address[myaddress]{Centre for Applicable Mathematics,Tata Institute of Fundamental Research, Post Bag No 6503, Sharadanagar, Bangalore - 560065, India.}
		\vspace{-2ex}
		\begin{abstract}
			We give a new  sufficient criteria to prove the uniqueness of the incompressible Euler equation in dimension $N\geq 2$. In their celebrated works by V. Scheffer \cite{Scheffer}, A. Shnirelman \cite{Shnirelman}, C. De Lellis and L. Székelyhidi Jr. \cite{Delellis}  they have obtained the nonuniqeness of weak solutions of incompressible Euler equation. Here we obtain uniqueness criteria for the same equation under some mild regularity condition on weak  solutions. Our proof is simple and can be employed to other equations like inhomogeneous incompressible Euler and Euler-Boussinesq equations. One of the key ingredients in our proof is commutator estimate \cite{CET,FGJ}. 
		\end{abstract}
		\begin{keyword}
			Incompressible Euler system \sep uniqueness \sep Besov space \sep energy conservation \sep commutator estimate
		\end{keyword}
		
	\end{frontmatter}
	\tableofcontents
	\section{Introduction}	
	In this article, we consider multidimensional incompressible Euler system, that is
	\begin{eqnarray}
	\pa_t\textbf{u}+\dv_x(\textbf{u}\otimes\textbf{u})+\nb_xp(x,t)&=&0\ \ \ \ \ \ \ \mbox{ for }(t,x)\in[0,T)\times\Omega,   \lab{I1}\\
	\dv_x\textbf{u}&=&0\ \ \ \ \ \ \ \mbox{ for }(t,x)\in[0,T)\times\Omega,    \lab{I3}\\
	\textbf{u}(0,x)&=&\textbf{u}_0(x)\hspace*{.2cm} \mbox{ for }x\in\Omega.\lab{ini}
	\end{eqnarray}
	Here $\Omega\subset\re^N$ is a bounded open set. In the system (\ref{I1})--(\ref{ini}), $\textbf{u}:[0,T)\times\Omega\rr\re^N$ represents velocity and the function $p:[0,T)\times\Omega\rr\re$ represents the pressure. This article concerns about a sufficient condition for uniqueness of  weak solutions to the system (\ref{I1})--(\ref{ini}).  We say, $\textbf{u}\in C([0,T),L^2(\Omega))$ is a weak solution to the system (\ref{I1})--(\ref{ini}) if it satisfies following integral identities:
	\begin{itemize}
		\item 
		\begin{equation}\lab{W1}
		\dip\domm\left[\textbf{u}\cdot\pa_t\psi+\textbf{u}\otimes\textbf{u}:\nb_x\psi\right]dxdt=\int\limits_{\Omega}\textbf{u}\cdot\psi(\tau_2,x)dx-\int\limits_{\Omega}\textbf{u}\cdot\psi(\tau_1,x)dx
		\end{equation}
		for $\psi\in C_c^{1}([0,T)\times\Omega,\re^N)$ with $\dv_x(\psi)=0$ and $0\leq\tau_1<\tau_2<T$.
		\item
		\begin{equation}\lab{W2}
		\domm\textbf{u}\cdot\nb_x\phi dxdt=0
		\end{equation}
		for $\phi\in C_c^{1}([0,T)\times\Omega)$ and $0\leq\tau_1<\tau_2<T$.
	\end{itemize}
	\begin{definition}[\textbf{admissible solution}]
		We say a weak solution is \textit{admissible} if $\textbf{u}\in L^{\f}([0,T),L^2(\Omega))$ satisfies the following inequality
		\begin{equation}\lab{W3}
		\int\limits_{\Omega}|\textbf{u}|^2(\tau_2,x)dx\leq \int\limits_{\Omega}|\textbf{u}|^2(\tau_1,x)dx
		\end{equation}
		for $0\leq\tau_1<\tau_2<T$.
	\end{definition}
	
	In this article, we consider the domain $\Omega=([-1,1]_{\pm1})^N$ for $N\geq 2$ to get rid of kinematic boundary terms. Note that in the weak formulation (\ref{W1}) the pressure term does not appear. Once the solution $\textbf{u}$ is known $p$ can be recovered upto a constant from the following Poisson equation
	\begin{equation}
	-\De p=\dv_x\dv_x(\textbf{u}\otimes\textbf{u}).
	\end{equation}
	The classical solutions to (\ref{I1})--(\ref{ini}) exist locally in time for smooth data and they are uniquely determined by the initial data.  Motivated by Kolmogrov's theory of turbulence and other physical aspects of incompressible Euler system, it is preferable to consider weak solutions of the system (\ref{I1})--(\ref{ini}). Unlike the classical solution, weak solutions are no more unique. In 1993, Scheffer \cite{Scheffer} first obtained the nonuniqueness results for incompressible Euler system showing the existence of a compact support solution to incompressible Euler system. In 1997, Shnirelmman \cite{Shnirelman} gave another approach to show the nonuniqueness result. In \cite{Delellis}, De Lellis and Székelyhidi gave a new construction for such solutions by convex integration techniques. It has to be noted that all the solutions constructed in \cite{Scheffer,Shnirelman,Delellis} are discontinuous. In this article, we assume  weak solution to be in Besov space $B^{\al}_{p,\f}$ with $\al>\frac{1}{3}$ and $p\geq 3$ (see (\ref{Besov}) for definition). Since $\al p>1$ they are continuous function on $\Omega$ for a.e. $t\in(0,T)$. Note that this is a `mild assumption' since it does not require to have full derivative.

	\par In 1949, Onsager \cite{Onsager} conjectured about the energy conservation based the H\"{o}lder exponent. He predicted that a weak solution conserves energy if it is H\"{o}lder continuous with exponent more than $\frac{1}{3}$ and conservation of energy fails if it is continuous with H\"{o}lder exponent less than $\frac{1}{3}$. In \cite{CET}, Constantin et al. showed the positive result of Onsager's conjecture. In that article, authors proved that equality holds in (\ref{W3}) if $\textbf{u}\in L^3([0,t); B^{\al}_{p,\f}(\Omega))$ for $\al>\frac{1}{3}$. They used a \textit{ commutator estimate.} Note that here $B^{\al}_{p}(\Omega)$ denotes the Besov space defined as follows
	\begin{equation}\lab{Besov}
	B^{\al}_{p,\f}(\Omega):=\left\{h\in L^{p}(\Omega);|h|_{B^{\al}_{p,\f}}:=\sup\limits_{\bar{\Omega_1}\subset\Omega}\sup\limits_{\xi\in\re^N,\xi+\Omega_1\subset\Omega}\frac{\|h(\cdot+\xi)-h(\cdot)\|_{L^p(\Omega)}}{|\xi|^\al}<\f\right\}
	\end{equation}
	for $p\geq1$ and $\al\in(0,1)$ (see \cite{Triebel} for more on these spaces).
	\par The other part of Onsager's conjecture has been concerned in \cite{Eyink} with some partial results. It has been gradually improved in recent studies \cite{BDIS,BDS,BDV} by Buckmaster, De Lellis, Székelyhidi, Vicol and finally settled by Isett \cite{Isett}. In those results, they constructed a H\"{o}lder continuous weak solution for each energy profile $e(t)\geq0$. These results also give examples in the favour of non-uniqueness of system (\ref{I1})--(\ref{ini}). In contrast, our next theorem is dedicated to the uniqueness result in this special class of solution which conserves the energy.  We state the following for the class of weak solutions considered in \cite{CET} with an extra condition (\ref{condition1}).

	\begin{theorem}\lab{theorem1}
		Let $\textbf{u},\textbf{v}$ be two weak solutions to the system (\ref{I1})--(\ref{I3}) corresponding to the initial data $\textbf{u}_0$. We assume that 
		\begin{equation*}
		\textbf{u},\textbf{v}\in L^3([0,T); B^{\al}_{p,\f}(\Omega))\ \mbox{ for }\al>\frac{1}{3},p\geq 3.
		\end{equation*}
		Suppose  there is a non-negative function $\mathcal{C}\in L^1([0,T))$ such that the following holds
		\begin{equation}\lab{condition1}
		\int\limits_{\Omega}{ \left[ - \zeta \cdot \textbf{v}(\tau, \cdot) (\zeta \cdot \nabla_x) \phi  + \mathcal{C}(\tau) |\zeta|^2 \phi \right] } dxdt \geq 0\ \mbox{ for }\zeta\in \re^N,\tau\in[0,T)
		\end{equation}
		for each $\phi\in C_c^{\f}(\Omega)$ with $\phi\geq 0$.
		Then $\textbf{u}\equiv\textbf{v}$ in $[0,T)\times\Omega$.
	\end{theorem}

	\par Weak-strong uniqueness is one available result towards the well-posedness theory of incompressible Euler system. This method requires the existence of one strong solution. In dimension two ,existence and uniqueness is known (see \cite{Lions1}) for solutions with $\dip\int\limits_{\Omega}|u(t,x)|^2dx=\int\limits_{\Omega}|u_0(x)|^2dx$ for each $t\in[0,T)$ and data in a special Hilbert space which is defined as follows
	\begin{equation*}
	H:=\{\textbf{u}_0\in L^2(\Omega); \dv_x\textbf{u}_0=0 \mbox{ in the sense of distribution}\}.
	\end{equation*} 
	In \cite{BrDS}, it has been shown that measure valued solution to the incompressible Euler system has to coincide with a solution to (\ref{I1})--(\ref{I3}), $\textbf{u}\in C([0,t],L^2(\re^N))$ satisfying the following condition
	\begin{equation}
	\int\limits_0^T\|(\nb_x\textbf{u})_{sym}\|_{L^{\f}}<\f
	\end{equation}
	where $(\nb_x\textbf{u})_{sym}$ is the symmetric part of the matrix $\nb_x\textbf{u}$. As we have mentioned before the classical solutions are valid for small time and very few class of initial data. In this article, we consider a weak solution in appropriate Besov space. We further impose a one-sided Lipschitz condition on that particular solution and show that other weak solutions have to coincide with the Besov solution if they come from the same initial data. This is the content of our next theorem.

	\begin{theorem}\lab{theorem2}
		Let $\textbf{u}$ be an admissible solution to the system (\ref{I1})--(\ref{I3}) with initial data $\textbf{u}_0$. Let $\textbf{v}$ be another pair of weak solution to the system (\ref{I1})--(\ref{I3}) with initial data $\textbf{u}_0$ such that
		\begin{equation}\lab{theo1-besov}
		\textbf{v}\in L^2([0,T);B^{\al}_{p,\f}(\Omega))\ \mbox{ for }\al>\frac{1}{2}, p\geq2.
		\end{equation}
		We assume that there is a non-negative function $\mathcal{C}\in L^1([0,T))$ such that
		\begin{equation}\lab{lipcon}
		\int\limits_{\Omega}{ \left[ - \zeta \cdot \textbf{v}(\tau, \cdot) (\zeta \cdot \nabla_x) \phi  + \mathcal{C}(\tau) |\zeta|^2 \phi \right] }dxdt \geq 0\ \mbox{ for }\zeta\in \re^N,\tau\in[0,T)
		\end{equation}
		holds for each $\phi\in C_c^{\f}(\Omega)$ with $\phi\geq 0$.
		Then $\textbf{u}\equiv\textbf{v}$ in $[0,T)\times\Omega$.
	\end{theorem}
	\par In this article we consider weak solutions with appropriate Besov regularity and a one-sided bound condition like (\ref{condition1}). Then we obtain that other weak solutions has to coincide with it if they come from the same initial data. This is not exactly weak-strong uniqueness since we do not have strong solution.  Therefore, in general they need not satisfy the system (\ref{I1})--(\ref{I3}) pointwise. That is why it does not follow from the classical weak-strong uniqueness. We mollify solutions and the equation as well and then pass the limit to get back everything in terms of solutions. These results cover a wider class of weak solutions of incompressible Euler system .
	\par As we mentioned before proofs are based on mollifying equation and then passing to the limit via commutator estimate \cite{Lions}. We have already mentioned the application of commutator estimate in the proof of Onsager's conjecture for system (\ref{I1})--(\ref{I3}) (see \cite{CET}). A similar type result has been obtained for compressible Euler system in \cite{FeGwGwWi}. The uniqueness result for dissipative solutions to isentropic compressible Euler system has been obtained recently in \cite{FGJ} for a wider class of weak solution. With the similar technique an uniqueness result has been proved in \cite{GJ} for a broader class of weak solution to complete Euler system. For a good survey on weak-strong uniqueness for Euler system we refer interested reader to \cite{Wiedemann}. We refer \cite{FeiNo,FeiNo1} for weak-strong uniqueness results in the context of Navier-Stokes equation. See \cite{Dafermos} for similar results in hyperbolic system.
	\par In later part of the article, we show the application of our proof in some incompressible systems namely, inhomogeneous incompressible Euler system and Euler–Boussinesq equations. Though the main commutator estimate is same for these equations we prefer to treat them in separate sections due to technical reasons since they don't follow directly  from the incompressible homogeneous case.
	
	\section{Proof of Theorem \ref{theorem1} and \ref{theorem2}}
	For a technical reason we first present the proof of Theorem \ref{theorem2} and then we prove Theorem \ref{theorem1} in subsection \ref{proof1}.	
	\subsection{Proof of Theorem \ref{theorem2}}
	For this section we define $E$ as follows
	\begin{equation}\lab{E}
	E(\textbf{u}\mid \textbf{v}):=\frac{1}{2}|\textbf{u}-\textbf{v}|^2.
	\end{equation}
	Let $\eta_\e $ be the standard mollifier sequence. Define $\textbf{v}_\e:=\textbf{v}*\eta_\e$. Now mollifying the system (\ref{I1})--(\ref{I3}) for $\textbf{v}$ we get
	\begin{eqnarray}
	\pa_t(\textbf{v})_\e+\dv_x(\textbf{v}\otimes\textbf{v})_\e+\nb_xp(x,t)_\e&=&0,\lab{M1}\\
	\dv_x(\textbf{v}_\e)&=&0.\lab{M2}
	\end{eqnarray}
	Next we put $\psi=\textbf{v}_\e$ in (\ref{W1}) and get
	\begin{equation}\lab{P1}
	\dip\domm\left[\textbf{u}\cdot\pa_t\textbf{v}_\e+\textbf{u}\otimes\textbf{u}:\nb_x\textbf{v}_\e\right]dxdt=\dip\int\limits_{\Omega}\textbf{u}\cdot\textbf{v}_\e(\tau_2,x)dx-\int\limits_{\Omega}\textbf{u}\cdot\textbf{v}_\e(\tau_1,x)dx.
	\end{equation}
	Putting $\phi=\frac{1}{2}|\textbf{v}_\e|^2$ in (\ref{W2}) we get
	\begin{equation}\lab{P2}
	\domm\textbf{u}\cdot\nb_x\textbf{v}_\e\cdot\textbf{v}_\e dxdt=0.
	\end{equation}
	By virtue of (\ref{W3}), (\ref{P1}) and fundamental theorem of calculus we get
	\begin{eqnarray}
	& &\dom E(\textbf{u}\mid \textbf{v}_\e)(\tau_2,x)dx-\dom E(\textbf{u} \mid \textbf{v}_\e)(\tau_1,x)dx\nonumber\\
	&=&\left[\dom\left(\frac{1}{2}|\textbf{u}|^2-\textbf{u}\cdot\textbf{v}_\e+\frac{1}{2}|\textbf{v}_\e|^2\right)dx\right]_{t=\tau_1}^{t=\tau_2}\nonumber\\
	&\leq&-\dip\domm\left[\textbf{u}\cdot\pa_t\textbf{v}_\e+\textbf{u}\otimes\textbf{u}:\nb_x\textbf{v}_\e-\textbf{v}_\e\cdot\pa_t\textbf{v}_\e\right]dxdt.\lab{EN1}
	\end{eqnarray}
	By employing (\ref{M2}) and (\ref{P2}) in (\ref{EN1}) we get
	\begin{eqnarray*}
		& &\dom E(\textbf{u}\mid \textbf{v}_\e)(\tau_2,x)dx-\dom E(\textbf{u} \mid \textbf{v}_\e)(\tau_1,x)dx\\
		&\leq&\domm\left[ (\textbf{v}_\e-\textbf{u})\cdot\pa_t\textbf{v}_\e-\textbf{u}\cdot\nb_x\textbf{v}_\e\cdot\textbf{u}+\textbf{u}\cdot\nb_x\textbf{v}_\e\cdot\textbf{v}_\e\right]dxdt.
	\end{eqnarray*}
	After a rearrangement of the terms, we have
	\begin{eqnarray}
	& &\dom E(\textbf{u}\mid \textbf{v}_\e)(\tau_2,x)dx-\dom E(\textbf{u} \mid \textbf{v}_\e)(\tau_1,x)dx\nonumber\\
	&\leq&\domm\left[ \pa_t\textbf{v}_\e+\textbf{v}_\e\cdot\nb_x\textbf{v}_\e\right]\cdot(\textbf{v}_\e-\textbf{u})dxdt+\domm(\textbf{u}-\textbf{v}_\e)\cdot\nb_x\textbf{v}_\e\cdot(\textbf{v}_\e-\textbf{u})dxdt.\lab{EN2}
	\end{eqnarray}
	Using the equation (\ref{M1}) in (\ref{EN2}) we get
	\begin{eqnarray}
	& &\dom E(\textbf{u}\mid \textbf{v}_\e)(\tau_2,x)dx-\dom E(\textbf{u} \mid \textbf{v}_\e)(\tau_1,x)dx\nonumber\\
	&\leq&\domm\left[ -\dv_x(\textbf{v}\otimes\textbf{v})_\e-\nb_xp_\e+\textbf{v}_\e\cdot\nb_x\textbf{v}_\e\right]\cdot(\textbf{v}_\e-\textbf{u})dxdt\nonumber\\
	&+&\domm(\textbf{u}-\textbf{v}_\e)\cdot\nb_x\textbf{v}_\e\cdot(\textbf{v}_\e-\textbf{u})dxdt\nonumber\\
	&=&R_1^{\e}+R_2^{\e},\lab{M3}
	\end{eqnarray}
	where $R_1^{\e},R_2^{\e}$ are defined as follows
	\begin{eqnarray*}
		R_1^{\e}&:=&\domm\left[ -\dv_x(\textbf{v}\otimes\textbf{v})_\e-\nb_xp_\e+\textbf{v}_\e\cdot\nb_x\textbf{v}_\e\right]\cdot(\textbf{v}_\e-\textbf{u})dxdt,\\
		R_2^{\e}&:=&\domm(\textbf{u}-\textbf{v}_\e)\cdot\nb_x\textbf{v}_\e\cdot(\textbf{v}_\e-\textbf{u})dxdt.
	\end{eqnarray*}
	Now we first analyze $	R_1^{\e}$ as follows
	\begin{eqnarray*}
		R_1^{\e}&=&\domm\left[ -\dv_x(\textbf{v}\otimes\textbf{v})_\e-\nb_xp_\e+\textbf{v}_\e\cdot\nb_x\textbf{v}_\e\right]\cdot(\textbf{v}_\e-\textbf{u})dxdt,\\
		&=&\domm\left[ \dv_x(\textbf{v}_\e\otimes\textbf{v}_\e)-\dv_x(\textbf{v}\otimes\textbf{v})_\e\right]\cdot(\textbf{v}_\e-\textbf{u})dxdt\\
		&+&\domm\left[- \dv_x(\textbf{v}_\e\otimes\textbf{v}_\e)+\textbf{v}_\e\cdot\nb_x\textbf{v}_\e\right]\cdot(\textbf{v}_\e-\textbf{u})dxdt-\domm\nb_xp_\e\cdot(\textbf{u}-\textbf{v}_\e)dxdt\\
		&=&\domm\left[ \dv_x(\textbf{v}_\e\otimes\textbf{v}_\e)-\dv_x(\textbf{v}\otimes\textbf{v})_\e\right]\cdot(\textbf{v}_\e-\textbf{u})dxdt\\
		&-&\domm\dv_x(\textbf{v}_\e)\textbf{v}_\e\cdot(\textbf{v}_\e-\textbf{u})dxdt-\domm\nb_xp_\e\cdot(\textbf{u}-\textbf{v}_\e)dxdt
	\end{eqnarray*}
	By employing (\ref{W2}) with $\phi=p_\e$ and (\ref{M2}) we get
	\begin{equation}\lab{R1}
	R_1^{\e}=\domm\left[ \dv_x(\textbf{v}_\e\otimes\textbf{v}_\e)-\dv_x(\textbf{v}\otimes\textbf{v})_\e\right]\cdot(\textbf{v}_\e-\textbf{u})dxdt.
	\end{equation}
	Next we estimate $R_2^{\e}$ with the help of (\ref{condition1}). Note that if we put $\phi=\eta_{\e}(x-y)$ and $\zeta=\textbf{u}(t,x)-\V(t,x)$ in (\ref{condition1}) we get
	\begin{equation*}
	(\textbf{u}-\textbf{v}_\e)\cdot\nb_x\textbf{v}_\e\cdot(\textbf{v}_\e-\textbf{u})\geq-\mathcal{C}(t)|\textbf{u}-\V|^2.
	\end{equation*}
	This yields
	\begin{equation}\lab{R2}
	R_2^\e\leq \domm\mathcal{C}(t)E(\textbf{u}|\textbf{v}_\e)(t,x)dxdt.
	\end{equation}
	Clubbing (\ref{M3}) with (\ref{R1}) and (\ref{R2}) we get
	\begin{equation}\lab{M4}
	\begin{array}{lll}
	& &\dip\dom\dip E(\textbf{u}\mid\textbf{v}_\e)(\tau_2,x)dx-\dom E(\textbf{u} \mid \textbf{v}_\e)(\tau_1,x)dx\\
	&\leq&\dip\domm\mathcal{C}(t)E(\textbf{u}|\textbf{v}_\e)(t,x)dxdt+\domm\left[ \dv_x(\textbf{v}_\e\otimes\textbf{v}_\e)-\dv_x(\textbf{v}\otimes\textbf{v})_\e\right]\cdot(\textbf{v}_\e-\textbf{u})dxdt.
	\end{array}
	\end{equation}
	\begin{lemma}[\textit{Commutator estimate}]\lab{comm}
		Let $D, D_1$ be two bounded domain in $\re^d$ such that $\bar{D}\subset D_1$. Let $f: {D}_1 \to \re^M$ be defined as $f=({f}_1,\cdots,{f}_M)$ such that $f_j\in B^{\al}_{q,\f}$ for each $j\in\{1,\cdots,M\}$, $\al\in(0,1)$ and $q\geq2$. Let $\eta_\e$ be a standard mollifier sequence with $\mbox{ supp}(\eta_{\e}) \subset \{ |z| < \e \}$. Let 
		$h : Q \to \re$ be a $C^2$ function  where $Q$ is a convex domain containing  the range of $f$.
		Then 
		\begin{equation}
		\left\| \nabla_z h( f_\e ) - \nabla_z  h(f)_\e \right\|_{L^{\frac{q}{2}} (D; \re^d) }
		\leq \e^{2\al-1} C|f|_{B^{\al}_{q,\f}}^2
		\end{equation} 
		for $\nabla_z = (\partial_{z_1}, \dots, \partial_{z_d})$ and the constant $C=\|h\|_{C^2}$.
	\end{lemma}
	We omit the proof of Lemma \ref{comm}. We refer \cite{FGJ,GJ} for the proof of a similar version of Lemma \ref{comm}. 
	\par Invoking Lemma \ref{comm} in (\ref{M4}) we get
	\begin{equation}\lab{P3}
	\begin{array}{lll}
	& &\dip\dom E(\textbf{u}\mid \textbf{v}_\e)(\tau_2,x)dx-\dom E(\textbf{u} \mid \textbf{v}_\e)(\tau_1,x)dx\\
	&\leq&\dip\domm\mathcal{C}(t)E(\textbf{u}|\textbf{v}_\e)(t,x)dxdt+\e^{2\al-1} C_1(|v|_{B^{\al}_{p,\f}}).
	\end{array}
	\end{equation}
	Now we are all set to pass the limit $\e\rr0$ in (\ref{P3}) and get
	\begin{equation}
	\dip\dom E(\textbf{u}\mid \textbf{v})(\tau_2,x)dx-\dom E(\textbf{u} \mid \textbf{v})(\tau_1,x)dx\leq\dip\domm\mathcal{C}(t)E(\textbf{u}|\textbf{v})(t,x)dxdt.
	\end{equation}
	Employing Gronwall's inequality and then passing the limit $\tau_1\rr0$ we get $\textbf{u}\equiv\textbf{v}$. This completes the proof of Theorem \ref{theorem2}.\qed
	\subsection{Proof of Theorem \ref{theorem1}}\lab{proof1}
	Suppose that $\textbf{v}_\e:=\textbf{v}*\eta_\e$ and $\textbf{u}_\e:=\textbf{u}*\eta_\e$.
	Mollifying the system (\ref{I1})--(\ref{I3}) for $\textbf{u}$ and $\textbf{v}$ we get
	\begin{eqnarray}
	\pa_t(\textbf{u})_\e+\dv_x(\textbf{u}\otimes\textbf{u})_\e+\nb_xp(x,t)_\e&=&0,\lab{M11}\\
	\dv_x(\textbf{u}_\e)&=&0,\lab{M12}\\
	\pa_t(\textbf{v})_\e+\dv_x(\textbf{v}\otimes\textbf{v})_\e+\nb_xp(x,t)_\e&=&0,\lab{M21}\\
	\dv_x(\textbf{v}_\e)&=&0.\lab{M22}
	\end{eqnarray}
	Integrating (\ref{M11}) against $\textbf{v}_\e$ we get
	\begin{equation*}
	\domm \left[\pa_t(\textbf{u})_\e+\dv_x(\textbf{u}\otimes\textbf{u})_\e+\nb_xp(x,t)_\e\right]\cdot\textbf{v}_\e dxdt=0.
	\end{equation*}
	After integrating by parts and applying (\ref{M22}) we obtain
	\begin{equation}\lab{P11}
	\dip\domm\left[\textbf{u}_\e\cdot\pa_t\textbf{v}_\e+(\textbf{u}\otimes\textbf{u})_\e:\nb_x\textbf{v}_\e\right]dxdt=\dip\int\limits_{\Omega}\textbf{u}_\e\cdot\textbf{v}_\e(\tau_2,x)dx-\int\limits_{\Omega}\textbf{u}_\e\cdot\textbf{v}_\e(\tau_1,x)dx.
	\end{equation}
	Integrating (\ref{M11}) against $\frac{1}{2}|\textbf{v}_\e|^2$, we get
	\begin{equation*}
	\domm\dv_x(\textbf{u}_\e)\frac{1}{2}|\textbf{v}_\e|^2=0.
	\end{equation*}
	Again integrating by parts we have
	\begin{equation}\lab{P12}
	\domm\textbf{u}_\e\cdot\nb_x\textbf{v}_\e\cdot\textbf{v}_\e dxdt=0.
	\end{equation}
	Let $E(\textbf{u}_\e\mid\V)$ be defined as in (\ref{E}). Now fundamental theorem of integral calculus and the equation (\ref{P11}) yield
	\begin{eqnarray*}
		& &\dom E(\textbf{u}_\e\mid \textbf{v}_\e)(\tau_2,x)dx-\dom E(\textbf{u}_\e \mid \textbf{v}_\e)(\tau_1,x)dx\\
		&=&\dip\domm\left[\textbf{u}_\e\cdot\pa_t\textbf{u}_\e-\textbf{u}_\e\cdot\pa_t\textbf{v}_\e-(\textbf{u}\otimes\textbf{u})_{\e}:\nb_x\textbf{v}_\e+\textbf{v}_\e\cdot\pa_t\textbf{v}_\e\right]dxdt.
	\end{eqnarray*}
	By employing (\ref{M11}), (\ref{M21}) and (\ref{P12}) we have
	\begin{eqnarray*}
		& &\dom E(\textbf{u}_\e\mid \textbf{v}_\e)(\tau_2,x)dx-\dom E(\textbf{u}_\e \mid \textbf{v}_\e)(\tau_1,x)dx\\
		&=&\dip\domm\left[-\textbf{u}_\e\cdot\dv_x(\textbf{u}\otimes\textbf{u})_\e+(\textbf{u}_\e-\V)\cdot\dv_x(\textbf{v}\otimes\textbf{v})_\e\right]dxdt\\
		&+&\domm \left[(\textbf{u}\otimes\textbf{u})_\e-(\textbf{u}_\e\otimes\textbf{u}_\e)\right]:\nb_x \textbf{v}_\e dxdt-
		\domm \textbf{u}_\e\cdot\nb_x\textbf{v}_\e\cdot\textbf{u}_\e dxdt.
	\end{eqnarray*}
	After a rearrangement of terms we get
	\begin{eqnarray}
	& &\dom E(\textbf{u}_\e\mid \textbf{v}_\e)(\tau_2,x)dx-\dom E(\textbf{u}_\e \mid \textbf{v}_\e)(\tau_1,x)dx\nonumber\\
	&=&\dip\domm\left[\textbf{u}_\e\cdot\left((\dv_x(\textbf{u}_\e\otimes\textbf{u}_\e)-\dv_x(\textbf{u}\otimes\textbf{u})_\e\right)+(\textbf{u}_\e-\V)\cdot\nb_x\textbf{v}_\e\cdot\textbf{v}_\e\right]dxdt\nonumber\\
	&-&\domm\left[\textbf{v}_\e\cdot\nb_x\textbf{v}_\e\cdot\textbf{u}_\e+\textbf{u}_\e\cdot\nb_x\textbf{u}_\e\cdot\textbf{u}_\e\right]dxdt\nonumber\\
	&+&\domm \left[(\textbf{u}\otimes\textbf{u})_\e-(\textbf{u}_\e\otimes\textbf{u}_\e)\right]:\nb_x \textbf{v}_\e dxdt
	-\domm (\textbf{u}_\e-\textbf{v}_\e)\cdot\nb_x\textbf{v}_\e\cdot\textbf{u}_\e.\lab{P13}
	\end{eqnarray}
	Note that (\ref{M12}) implies
	\begin{equation}\lab{P4}
	\domm \nb_x\left(|\textbf{u}_\e|^2+|\textbf{v}_\e|^2\right)\cdot\textbf{u}_\e dxdt=0.
	\end{equation}
	Therefore applying (\ref{P4}) in (\ref{P3}), we have 
	\begin{eqnarray*}
		& &\dom E(\textbf{u}_\e\mid \textbf{v}_\e)(\tau_2,x)dx-\dom E(\textbf{u}_\e \mid \textbf{v}_\e)(\tau_1,x)dx\\
		&=&\dip\domm\left[\textbf{u}_\e\cdot\left(\dv_x(\textbf{u}_\e\otimes\textbf{u}_\e)-\dv_x(\textbf{u}\otimes\textbf{u})_\e\right)\right]dxdt\\
		&-&\domm(\textbf{u}_\e-\V)\cdot\nb_x\textbf{v}_\e\cdot(\textbf{u}_\e-\V) dx dt\\
		&+&\domm \left[(\textbf{u}\otimes\textbf{u})_\e-(\textbf{u}_\e\otimes\textbf{u}_\e)\right]:\nb_x \textbf{v}_\e dxdt.
	\end{eqnarray*}
	Again integrating by parts we get
	\begin{eqnarray*}
		& &\dom E(\textbf{u}_\e\mid \textbf{v}_\e)(\tau_2,x)dx-\dom E(\textbf{u}_\e \mid \textbf{v}_\e)(\tau_1,x)dx\\
		&=&\domm \left[(\textbf{u}\otimes\textbf{u})_\e-(\textbf{u}_\e\otimes\textbf{u}_\e)\right]:\nb_x (\textbf{v}_\e-\textbf{u}_\e) dxdt\\
		&-&\domm(\textbf{u}_\e-\V)\cdot\nb_x\textbf{v}_\e\cdot(\textbf{u}_\e-\V) dx dt.
	\end{eqnarray*}
	By virtue of (\ref{condition1}) we have 
	\begin{equation}
	\dom E(\textbf{u}_\e\mid \textbf{v}_\e)(\tau_2,x)dx-\dom E(\textbf{u}_\e \mid \textbf{v}_\e)(\tau_1,x)dx\leq\domm\mathcal{C}(t) E(\textbf{u}_\e\mid \textbf{v}_\e)(t,x)dxdt+\mathcal{R}_\e
	\end{equation}
	where $\mathcal{R}_\e$ is defined as follows
	\begin{equation}\lab{RE}
	\mathcal{R}_\e:=\domm \left[(\textbf{u}\otimes\textbf{u})_\e-(\textbf{u}_\e\otimes\textbf{u}_\e)\right]:\nb_x (\textbf{v}_\e-\textbf{u}_\e) dxdt.
	\end{equation}
	\begin{lemma}[Constantin et al. \cite{CET}]\lab{commutator1} Let $\Omega\subset\re^N$ be a bounded domain and $T>0$. Let $\textbf{u},\textbf{v}\in B^{\al}_{p,\f}(\Omega,\re^N)$ with $\al\in(0,1)$ and $p\geq 3$. Then the following estimate holds
		\begin{equation*}
		\left|\domm \left[(\textbf{u}\otimes\textbf{u})_\e-(\textbf{u}_\e\otimes\textbf{u}_\e)\right]:\nb_x (\textbf{v}_\e-\textbf{u}_\e) dxdt\right|\leq\e^{3\al-1}C(\Omega,T)|\textbf{u}|_{B^{\al}_{p,\f}}^2(|\textbf{u}|_{B^{\al}_{p,\f}}+|\textbf{v}|_{B^{\al}_{p,\f}})
		\end{equation*}
		for $0<\tau_1<\tau_2<T$.
	\end{lemma}
	Proof of this lemma can be found in \cite[page 208-209]{CET}. 
	\par By applying Lemma \ref{commutator1} in (\ref{RE}) we get the following estimate of $\mathcal{R}_{\e}$ 
	\begin{equation*}
	|\mathcal{R}_\e|\leq \e^{3\al-1}C(\Omega,T)|\textbf{u}|_{B^{\al}_{p,\f}}^2(|\textbf{u}|_{B^{\al}_{p,\f}}+|\textbf{v}|_{B^{\al}_{p,\f}}).
	\end{equation*}
	This yields
	\begin{eqnarray*}
		\dom E(\textbf{u}_\e\mid \textbf{v}_\e)(\tau_2,x)dx-\dom E(\textbf{u}_\e \mid \textbf{v}_\e)(\tau_1,x)dx&\leq&\domm\mathcal{C}(t) E(\textbf{u}_\e\mid \textbf{v}_\e)(t,x)dxdt\\
		&+&\e^{3\al-1}C(\Omega,T)|\textbf{u}|_{B^{\al}_{p,\f}}^2(|\textbf{u}|_{B^{\al}_{p,\f}}+|\textbf{v}|_{B^{\al}_{p,\f}}).
	\end{eqnarray*}
	Passing the limit $\e\rr0$ we get
	\begin{eqnarray*}
		\dom E(\textbf{u}\mid \textbf{v})(\tau_2,x)dx-\dom E(\textbf{u} \mid \textbf{v})(\tau_1,x)dx&\leq&\domm\mathcal{C}(t) E(\textbf{u}\mid \textbf{v})(t,x)dxdt.
	\end{eqnarray*}
	Now passing the limit $\tau_1\rr0$ and invoking Gronwall's inequality we get $\textbf{u}\equiv\textbf{v}$. This completes the proof of Theorem \ref{theorem1}.\qed
	
	\section{Inhomogeneuos incompressible Euler system}
	In this section, we consider the inhomogeneous incompressible Euler system 
	\begin{eqnarray}
	\pa_t(\vr\textbf{u})+\dv_x(\vr\textbf{u}\otimes\textbf{u})+\nb_xp(x,t)&=&0 \hspace{2.2cm} \mbox{ for }(t,x)\in[0,T)\times\Omega,\lab{I31}\\
	\pa_t\vr+\dv_x(\vr\textbf{u})&=&0 \hspace{2.2cm} \mbox{ for }(t,x)\in[0,T)\times\Omega,\lab{I32}\\
	\dv_x(\textbf{u})&=&0 \hspace{2.2cm} \mbox{ for }(t,x)\in[0,T)\times\Omega,\lab{I33}\\
	(\vr,\textbf{u})(0,x)&=&(\vr_0,\textbf{u}_0)(x)\hspace{0.5cm}\mbox{ for }x\in\Omega.
	\end{eqnarray}
	First we define the weak solution to the system (\ref{I31})--(\ref{I33}) in a similar way as we have done for the system (\ref{I1})--(\ref{I3}).
	
	\noindent\textbf{Weak formulation:} We say $(\rho,\textbf{u})\in C([0,T),L^2(\Omega))$ is a weak solution to (\ref{I31})--(\ref{I33}) if it satisfies following integral equations
	\begin{itemize}
		\item 
		\begin{equation}\lab{W31}
		\dip\domm\left[\vr\textbf{u}\cdot\pa_t\psi+\vr\textbf{u}\otimes\textbf{u}:\nb_x\psi\right]dxdt\\
		=\dip\int\limits_{\Omega}\vr\textbf{u}\cdot\psi(\tau_2,x)dx-\int\limits_{\Omega}\vr\textbf{u}\cdot\psi(\tau_1,x)dx
		\end{equation}
		for $\psi\in C_c^{1}([0,T)\times\Omega,\re^N)$ and $0\leq\tau_1<\tau_2<T$ with $\dv_x(\psi)=0$.
		\item 
		\begin{equation}\lab{W32}
		\domm\left[\vr\pa_t\phi+\vr\textbf{u}\cdot\nb_x\phi\right]=\int\limits_{\Omega}\vr\phi(\tau_2,x)dx-\int\limits_{\Omega}\vr\phi(\tau_1,x)dx
		\end{equation}
		for $\phi\in C_c^{1}([0,T)\times\Omega)$ nd $0\leq\tau_1<\tau_2<T$.	
		\item
		\begin{equation}\lab{W33}
		\domm\textbf{u}\cdot\nb_x\phi dxdt=0
		\end{equation}
		for $\phi\in C_c^{1}([0,T)\times\Omega)$ and $0\leq\tau_1<\tau_2<T$.
	\end{itemize}
	\begin{definition}[\textbf{admissible solution}]
		We say a weak solution to the system (\ref{I31})--(\ref{I33}) is admissible if it satisfies the following inequality
	\end{definition}
	\begin{equation}
	\int\limits_{\Omega}\vr|\textbf{u}|^2(\tau_2,x)dx\leq \int\limits_{\Omega}\vr|\textbf{u}|^2(\tau_1,x)dx.
	\end{equation}
	
	\begin{theorem}\lab{theorem3}
		Let $(\vr,\textbf{u}),(r,\textbf{v})$ be two weak solutions to the system (\ref{I31})--(\ref{I33}) with same initial data $(\vr_0,\textbf{u}_0)$  such that the following holds 
		\begin{equation*}
		\vr, r,\vr\textbf{u},r\textbf{v}, \textbf{u},\textbf{v} \in B^{\al}_{p,\f}([0,T)\times\Omega), \mbox{ with }\al>\frac{1}{3},p\geq 3.
		\end{equation*}
		Suppose  there is a non-negative function $\mathcal{C}\in L^1([0,T))$ such that the following holds
		\begin{equation}\lab{condition2}
		\int\limits_{\Omega}{ \left[ - \zeta \cdot \textbf{v}(\tau, \cdot) (\zeta \cdot \nabla_x) \phi  + \mathcal{C}(\tau) |\zeta|^2 \phi \right] } dxdt \geq 0\ \mbox{ for }\zeta\in \re^N,\tau\in[0,T)
		\end{equation}
		for each $\phi\in C_c^{\f}(\Omega)$ with $\phi\geq 0$.
		Then we have $\vr\equiv r$ and $\textbf{u}\equiv\textbf{v}$  in $[0,T)\times\Omega$.
	\end{theorem}
	\subsection{Proof of Theorem \ref{theorem3}}
	In the context of inhomogeneous incompressible Euler system we define 
	\begin{equation*}
	E_1(\textbf{u}\mid \textbf{v}):=\frac{1}{2}\vr|\textbf{u}-\textbf{v}|^2.
	\end{equation*}
	By a similar method as we have done in the proof of Theorem \ref{theorem2}, we can prove $\textbf{u}\equiv\textbf{v}$. Next we mollify (\ref{I32}) and (\ref{I33}) for $(\vr,\textbf{u})$ and $(r,\textbf{u})$ to get
	\begin{eqnarray}
	\pa_t\vr_\e+\dv_x(\vr\textbf{u})_\e&=&0,\lab{M31}\\
	\dv_x(\textbf{u}_\e)&=&0,\lab{32}\\
	\pa_tr_\e+\dv_x(r\textbf{u})_\e&=&0.\lab{M33}
	\end{eqnarray}
	Subtracting (\ref{M33}) from (\ref{M31}) and multiplying by $(\vr_\e-r_\e)$ we get
	\begin{eqnarray*}
		\pa_t\frac{1}{2}|\vr_\e-r_\e|^2+\dv_x\left((\vr\textbf{u})_\e-(r\textbf{u})_\e\right)(\vr_\e-r_\e)=0.
	\end{eqnarray*}
	This yields
	\begin{eqnarray*}
		& &\dom\frac{1}{2}|\vr_\e-r_\e|^2(\tau_2,x)dx-\dom\frac{1}{2}|\vr_\e-r_\e|^2(\tau_2,x)dx\\
		&=&-\domm\dv_x\left((\vr\textbf{u})_\e-(r\textbf{u})_\e\right)(\vr_\e-r_\e)dxdt.
	\end{eqnarray*}
	Integrating by parts we get
	\begin{eqnarray*}
		& &\dom\frac{1}{2}|\vr_\e-r_\e|^2(\tau_2,x)dx-\dom\frac{1}{2}|\vr_\e-r_\e|^2(\tau_1,x)dx\\
		&=&\domm\left((\vr\textbf{u})_\e-(r\textbf{u})_\e\right)\cdot\nb_x(\vr_\e-r_\e)dxdt.
	\end{eqnarray*}
	After a rearrangement of terms we have
	\begin{eqnarray}
	& &\dom\frac{1}{2}|\vr_\e-r_\e|^2(\tau_2,x)dx-\dom\frac{1}{2}|\vr_\e-r_\e|^2(\tau_1,x)dx\nonumber\\
	&=&\domm\left((\vr\textbf{u})_\e-(\vr_\e\textbf{u}_\e)+(r_\e\textbf{u}_\e)-(r\textbf{u})_\e\right)\cdot\nb_x(\vr_\e-r_\e)dxdt\nonumber\\
	&+&\domm\left((\vr_\e\textbf{u}_\e)-(r_\e\textbf{u}_\e)\right)\cdot\nb_x(\vr_\e-r_\e)dxdt\nonumber\\
	&=&\domm\left((\vr\textbf{u})_\e-(\vr_\e\textbf{u}_\e)+(r_\e\textbf{u}_\e)-(r\textbf{u})_\e\right)\cdot\nb_x(\vr_\e-r_\e)dxdt\nonumber\\
	&+&\domm\textbf{u}_\e\cdot\nb_x(\frac{1}{2}|\vr_\e-r_\e|^2)dxdt.\lab{P44}
	\end{eqnarray}
	Note that the last term in (\ref{P44}) vanishes by virtue of (\ref{W33})  with $\phi=\frac{1}{2}|\vr_\e-r_\e|^2$. Therefore we get
	\begin{eqnarray*}
		& &\dom\frac{1}{2}|\vr_\e-r_\e|^2(\tau_2,x)dx-\dom\frac{1}{2}|\vr_\e-r_\e|^2(\tau_1,x)dx\\
		&=&\domm\left((\vr\textbf{u})_\e-(\vr_\e\textbf{u}_\e)+(r_\e\textbf{u}_\e)-(r\textbf{u})_\e\right)\cdot\nb_x(\vr_\e-r_\e)dxdt.
	\end{eqnarray*}
	By using Lemma \ref{comm} we pass the limit $\e\rr0$ to get
	\begin{equation}
	\dom\frac{1}{2}|\vr_\e-r_\e|^2(\tau_2,x)dx=\dom\frac{1}{2}|\vr_\e-r_\e|^2(\tau_1,x)dx.
	\end{equation}
	Next we pass the limit $\tau_1\rr 0$. Therefore we have $\vr\equiv r$.
	This completes the proof of Theorem \ref{theorem3}.\qed
	\section{Euler–Boussinesq equations}
	This section deals with an uniqueness result for \textit{Euler–Boussinesq equations} which is the following
	\begin{eqnarray}
	\pa_t\textbf{u}+\dv_x(\textbf{u}\otimes\textbf{u})+\nb_xp(x,t)&=&\theta g\hspace{2cm} \mbox{ for }(t,x)\in[0,T)\times\Omega,\lab{I41}\\
	\pa_t\theta+\dv_x(\theta\textbf{u})&=&0\hspace{2.2cm} \mbox{ for }(t,x)\in[0,T)\times\Omega,\lab{I42}\\
	\dv_x(\textbf{u})&=&0\hspace{2.2cm} \mbox{ for }(t,x)\in[0,T)\times\Omega,\lab{I43}\\
	(\theta,\textbf{u})(0,x)&=&(\theta_0,\textbf{u}_0)(x)\hspace{0.55cm}\mbox{ for }x\in\Omega.
	\end{eqnarray}
	We say $(\theta,\textbf{u})\in C([0,T),L^2(\Omega))$ is a weak solution to the system (\ref{I41})--(\ref{I43}) if it satisfies the following integral equations
	\begin{itemize}
		\item 
		\begin{equation}\lab{W41}
		\dip\domm\left[\textbf{u}\cdot\pa_t\psi+\textbf{u}\otimes\textbf{u}:\nb_x\psi+\theta g\cdot\textbf{u}\right]dxdt\\
		=\dip\int\limits_{\Omega}\textbf{u}\cdot\psi(\tau_2,x)dx-\int\limits_{\Omega}\textbf{u}\cdot\psi(\tau_1,x)dx
		\end{equation}
		for $\psi\in C_c^{1}([0,T)\times\Omega,\re^N)$ and $0\leq\tau_1<\tau_2<T$ with $\dv_x(\psi)=0$.
		\item 
		\begin{equation}\lab{W42}
		\domm\left[\theta\pa_t\phi+\theta\textbf{u}\cdot\nb_x\phi\right]=\int\limits_{\Omega}\theta\phi(\tau_2,x)dx-\int\limits_{\Omega}\theta\phi(\tau_1,x)dx
		\end{equation}
		for $\phi\in C_c^{1}([0,T)\times\Omega)$ nd $0\leq\tau_1<\tau_2<T$.	
		\item
		\begin{equation}
		\domm\textbf{u}\cdot\nb_x\phi dxdt=0
		\end{equation}
		for $\phi\in C_c^{1}([0,T)\times\Omega)$ and $0\leq\tau_1<\tau_2<T$.
	\end{itemize}
	
	By a similar argument as we have given in proofs of Theorem \ref{theorem1} and \ref{theorem3}, we can prove the following theorem.
	\begin{theorem}\lab{theorem4}
		Let $(\theta_1,\textbf{u}_1),(\theta_2,\textbf{u}_2)$ be two weak solution to the system (\ref{I41})--(\ref{I43}) such that the following holds
		\begin{equation}
		\theta_1,\theta_2,\theta_1\textbf{u}_1,\theta_2\textbf{u}_2,\textbf{u}_1,\textbf{u}_2\in B^{\al}_{p,\f}((0,T)\times\Omega), \mbox{ with }\al>\frac{1}{3},p\geq 3.
		\end{equation}
		Additionally, we assume that there is a non-negative function $\mathcal{C}\in L^1([0,T))$ such that the following holds
		\begin{equation}\lab{condition2}
		\int\limits_{\Omega}{ \left[ - \zeta \cdot \textbf{u}_2(\tau, \cdot) (\zeta \cdot \nabla_x) \phi  + \mathcal{C}(\tau) |\zeta|^2 \phi \right] } dxdt \geq 0\ \mbox{ for }\zeta\in \re^N,\tau\in[0,T)
		\end{equation}
		for each $\phi\in C_c^{\f}(\Omega)$ with $\phi\geq 0$.
		Then we have $\theta_1\equiv \theta_2$ and $\textbf{u}_1\equiv\textbf{u}_2$  in $[0,T)\times\Omega$.
		
	\end{theorem}

	\noindent\textbf{Acknowledgement.}  The first author would like to thank Inspire faculty-research grant\\ DST/INSPIRE/04/2016/000237.

\end{document}